\documentclass[a4paper,11pt,reqno,twoside]{amsart}

\usepackage{amsmath}
\usepackage{amsfonts}
\usepackage{amssymb}
\usepackage{amsthm}
\usepackage{color}
\usepackage{ifpdf}
\usepackage{array}
\usepackage{url}
\usepackage{multirow,bigdelim}
\usepackage{ascmac}

\numberwithin{equation}{section}

\newtheorem{theorem}{Theorem}[section]

\newtheorem{corollary}{Corollary}[section]

\newcommand*{\C}{\mathbb{C}}
\newcommand*{\R}{\mathbb{R}}

\newcommand{\comment}[1]{}
\title[Screw functions of Dirichlet series]%
      {Screw functions of Dirichlet series \\ in the extended Selberg class} 
\author[M. Suzuki]{Masatoshi Suzuki}
\date{Version of \today}
\AtBeginDocument{%
\begin{abstract}
We introduce screw functions for Dirichlet series in the extended Selberg class. 
We then prove that the Grand Riemann Hypothesis 
for a member of the extended Selberg class 
is equivalent to the non-positivity of the corresponding screw function. 
\end{abstract}
\maketitle
}
\begin{document}

\thispagestyle{headings}

%
\section{Introduction} 
%

In \cite{Su22}, we introduced a screw function corresponding to the Riemann zeta-function 
whose pointwise non-positivity is equivalent to the Riemann hypothesis 
for the Riemann zeta-function. 
As a large class of functions similar to the Riemann zeta-function, 
the Selberg class $\mathcal{S}$ and 
the extended Selberg class $\mathcal{S}^{\sharp}$ 
are well-known and important in number theory. 
The Riemann zeta-function, 
Dirichlet $L$-functions, 
and $L$-functions of holomorphic cusp forms 
are typical members of $\mathcal{S}$ (cf. \cite{Pe05}). 
%
Therefore, it is naturally expected that 
a member $F$ of  the extended Selberg class  
also corresponds to a screw function similarly to the Riemann zeta-function. 
This paper implements that expectation.

Let $F$ be a member of the extended Selberg class $\mathcal{S}^{\sharp}$. 
The nontrivial zeros of $F(s)$ mean the zeros of the entire function $\xi_F(s)$, 
which is the Selberg class analogue of the Riemann xi-function 
defined in \eqref{s203} below. 
We call $\gamma \in \C$ a nontrivial zero of $F(1/2-iz)$ 
if $1/2-i\gamma$ is a nontrivial zero of $F(s)$. 
The Grand Riemann Hypothesis (GRH, for short) for $F$ 
is the assertion that all nontrivial zeros of $F(s)$ lie on the line $\Re(s)=1/2$, 
which is equivalent to the statement that all nontrivial zeros of $F(1/2-iz)$ are real. 
We define the function $g_F(t)$ on the real line by 
\begin{equation} \label{eq_101} 
g_F(t):= 
- iB_F t
- \frac{m_0}{2}\,t^2 + \sum_{\gamma \in \Gamma_F\setminus\{0\}} 
m_\gamma \frac{e^{-i\gamma t}-1}{\gamma^2},  
\end{equation}
%
where $\Gamma_F$ denotes the set of all nontrivial zeros of $F(1/2-iz)$, 
$m_\gamma$ represents the multiplicity of $\gamma \in \Gamma_F$, 
$m_0=0$ if $0 \not \in \Gamma_F$, and $B_F$ is the constant in \eqref{eq_0123_1} below.
For ease of introduction, 
definition \eqref{eq_101} refers to 
\cite[(1.3)]{Su22} 
rather than definition \cite[(1.1)]{Su22} of $\Psi(t)=-g(t)$, 
because 
we need the 
semi-extended Selberg class $\mathcal{S}^{\sharp\flat}\,(\subset \mathcal{S}^{\sharp})$ 
introduced by Smajlovi\'{c} \cite{Sm10} (and reviewed in Section \ref{section_2}) 
to define $g_F(t)$ without using zeros (cf. Theorem \ref{thm_401} below). 
The difference from the case of the Riemann zeta-function 
is the existence of the term $-m_0t^2/2$ corresponding to $\gamma=0$ 
and the lack of symmetry $\gamma \mapsto -\overline{\gamma}$ 
in the sum over zeros.
We call $g_F(t)$ the screw function of $F$ 
because 
it is a screw function in the sense of \cite{KrLa14} 
assuming the GRH for $F$ is true (see Section \ref{section_5} below). 
%
The main result of the present paper is the following generalization of \cite[Theorem 1.7]{Su22}.

\begin{theorem} \label{thm_101}
Let $F$ be a member of 
$\mathcal{S}^{\sharp}$ 
and let $g_F(t)$ be the screw function defined by \eqref{eq_101}. 
We assume that $F(s)$ has no real zeros except for the possible zero at $s=1/2$. 
Then, the GRH for $F$ is true if and only if 
$\Re(-g_F(t)) \geq 0$ for all $t \geq t_0$ for some $t_0 \geq 0$. 
\end{theorem} 


\begin{corollary} \label{cor_1}
Let $F \in \mathcal{S}$. Suppose that $\Re(-g_F(t)) \geq 0$ for all $t \geq t_0$ 
and $s=1$ is not a zero of $F(s)$. 
Then there exists $\delta>0$ such that 
$F(s)$ has no zeros in the right-half plane $\Re(s)>1-\delta$. 
\end{corollary}

In the following, 
we review the Selberg class in Section \ref{section_2}, 
and prove Theorem \ref{thm_101} in Section \ref{section_3}. 
Then we show that the screw function $g_F$ allows a representation 
without zeros if we restrict $F$ to 
the semi-extended Selberg class $\mathcal{S}^{\sharp\flat}$ in Section 4.
Finally, we comment mainly on generalizations 
of other results in \cite{Su22} in Section \ref{section_5}. 

\section{The Selberg class and explicit formulas} \label{section_2}

\subsection{} 

The Selberg class $\mathcal{S}$ introduced by A. Selberg in 1992 
consists of the Dirichlet series 
\begin{equation} \label{s201}
F(s) = \sum_{n=1}^{\infty} \frac{a_F(n)}{n^s}
\end{equation}
satisfying the following five axioms: 
\begin{enumerate}
\item[(S1)] The Dirichlet series \eqref{s201} converges absolutely if $\Re(s)>1$. 
\item[(S2)] There exists an integer $m \geq 0$ 
such that $(s - 1)^mF(s)$ extends to an entire function of finite order. 
The smallest $m$ is denoted by $\mathfrak{m}_F$. 
\item[(S3)] $F$ satisfies the functional equation
\begin{equation} \label{s202}
\xi_F(s)=\omega \, \overline{\xi_F(1-\bar{s})},
\end{equation}
where
\begin{equation} \label{s203}
\aligned 
\xi_F(s) 
& := s^{\mathfrak{m}_F}(s-1)^{\mathfrak{m}_F}Q^s \prod_{j=1}^{r} \Gamma(\lambda_{j} s + \mu_{j}) \, F(s) \\
& = s^{\mathfrak{m}_F}(s-1)^{\mathfrak{m}_F}\gamma(s)F(s), 
\endaligned 
\end{equation}
$\Gamma(s)$ is the gamma function,  
and $r \geq 0$, $Q > 0$, $\lambda_{j} > 0$, $\mu_{j} \in \C$ with $\Re(\mu_{j})\geq 0$, 
$\omega \in \C$ with  $|\omega| = 1$ are parameters depending on $F$. 
\item[(S4)] For every $\varepsilon > 0$, $a_F(n) \ll_\varepsilon n^\varepsilon$. 
\item[(S5)] 
$\log F(s) = \sum_{n=1}^{\infty} b_F(n)\,n^{-s}$,  
where $b_F(n) = 0$ unless $n = p^k$ with $k \geq 1$, and $b_F(n) \ll n^\theta$ for some $\theta < 1/2$.
\end{enumerate}

From (S3) and (S5), 
$F \in \mathcal{S}$ has no zeros outside the critical strip $0 \leq \Re(s) \leq 1$ 
except for zeros in the left half-plane $\Re(s)  \leq 0$ located at poles of the involved gamma factors 
(\cite[p. 29]{Pe05}). 
The zeros lying in the critical strip are called the nontrivial zeros. 
The nontrivial zeros are infinitely many unless $F \equiv 1$ 
and coincide with the zeros of the entire function $\xi_F(s)$ of \eqref{s203}. 
The extended Selberg class $\mathcal{S}^\sharp$ 
is the class of functions satisfying (S1)--(S3) above. 
Note that the data $(\omega, Q,r,\lambda_j, \mu_j)$ 
in (S3) are not uniquely determined by $F \in \mathcal{S}^\sharp$.  

For $F \in \mathcal{S}^\sharp$, we define 
$F^\ast(s):=\overline{F(\bar{s})}$. 
Then, $F^\ast \in \mathcal{S}^\sharp$ with 
$a_{F^\ast}(n)=\overline{a_F(n)}$, 
$\mathfrak{m}_{F^\ast}=\mathfrak{m}_F$, 
$\omega_{F^\ast}(n)=\overline{\omega_F(n)}$, 
$r_{F^\ast}=r_F$, 
$Q_{F^\ast}=Q_F$, 
$\lambda_{F^\ast,j}=\lambda_{F,j}$, and 
$\mu_{F^\ast,j}=\overline{\mu_{F,j}}$, 
where the chosen $\omega$, $Q$, $r$, $\lambda_j$, and $\mu_j$ in (S3) 
are written as $\omega_F$, $Q_F$, $r_F$, $\lambda_{F,j}$, and $\mu_{F,j}$, respectively. 
In particular, 
\begin{equation} \label{s204}
\xi_{F^\ast}(s) = \overline{\xi_F(\bar{s})}. 
\end{equation}
The survey \cite{Pe05} and its sequel are good introductions to 
the theory of the Selberg class and the extended Selberg class.

To represent the screw function $g_F(t)$ without using zeros of $F(1/2-iz)$, 
we recall the subclass $\mathcal{S}^{\sharp\flat}$ 
of the extended Selberg class $\mathcal{S}^{\sharp}$ 
introduced by Smajlovi\'{c} in \cite[Section 2]{Sm10} 
that consists of $F \in \mathcal{S}^{\sharp}$ 
satisfying the Euler sum condition: 
\begin{enumerate}
\item[(S5')] The logarithmic derivative of $F$ possesses a Dirichlet series representation
\[
-\frac{F'}{F}(s)=\sum_{n=2}^{\infty} \frac{c_F(n)}{n^s}
\]
converging absolutely for $\Re(s)>1$. 
\end{enumerate}
The subclass $\mathcal{S}^{\sharp\flat}$, which we refer to as the {\it semi-extended Selberg class}, 
contains the Selberg class $\mathcal{S}$ (\cite[Theorem 2.1]{Sm10}).
The Euler sum condition (S5') is satisfied 
if $\log F(s)$ admits a Dirichlet series with suitable convergence. 
For conditions related to the Dirichlet series representation of $\log F(s)$, 
see \cite[Theorem 11.14]{Ap76} for example, 
and also \cite[Satz 12]{La33} regarding the convergence of the Dirichlet series. 

\subsection{} \label{section_2_2}

For $F \in \mathcal{S}^{\sharp}$, 
the function 
$\xi_F(s)$ is an entire function of order one (\cite[Lemma 3.3 and a note below the proof]{Sm10}). 
Therefore, the Hadamard factorization theorem implies that 
$\xi_F(1/2-iz)$ has the product  formula 
\[
\xi_F\left(\frac{1}{2}-iz \right)
= e^{A_F+B_Fz} z^{m_0} \prod_{\gamma \in \Gamma_F\setminus\{0\}}
\left( 1 - \frac{z}{\gamma}
\right)^{m_\gamma} \exp\left(\frac{m_\gamma z}{\gamma}\right). 
\]
Taking the logarithmic derivative, 
\begin{equation} \label{eq_0914_1}
\frac{\xi_F^\prime}{\xi_F}\left( \frac{1}{2}-iz\right) 
= i B_F +\frac{im_0}{z} 
+ i \sum_{\gamma \in \Gamma_F\setminus\{0\}} 
m_\gamma \left( \frac{1}{z-\gamma} + \frac{1}{\gamma} \right),
\end{equation}
where the sum on the right-hand side converges absolutely and uniformly 
on every compact subset of $\C\setminus \Gamma_F$. 
By the functional equation \eqref{s202} and \eqref{eq_0914_1}, 
we find that $B_F$ is real, and 
\begin{equation} \label{eq_0123_1}
iB_F = 
 \lim_{z \to 0}\left[\frac{\xi_F^\prime}{\xi_F}\left(\frac{1}{2}-iz\right) - \frac{im_0}{z}\right].
\end{equation}

\section{Proof of the main results} \label{section_3}

\subsection{Proof of Theorem \ref{thm_101}} 

Assuming the GRH for $F$ is true,  we have 
\[
\Re(-g_F(t)) 
= \frac{m_0}{2}\,t^2 + 
\sum_{\gamma \in \Gamma_F\setminus\{0\}} m_\gamma \frac{1-\cos(\gamma t)}{\gamma^2}, 
\]
since $iB_F$ is pure imaginary as observed in Section \ref{section_2_2}. 
The values of $\Re(-g_F(t))$ are clearly non-negative as seen from the right-hand side. 

Now we prove the converse claim. 
For each term on the right-hand side of \eqref{eq_0914_1}, we have 
\[
\int_{0}^{\infty} \frac{e^{-i\gamma t}-1}{\gamma^2}\, e^{izt} \, dt
= \frac{i}{z^2}\left( \frac{1}{z-\gamma} + \frac{1}{\gamma}\right) 
\quad \text{if}~\Im(z) > \Im(\gamma),  
\]
and 
\begin{equation} \label{eq_0913_5} 
\int_{0}^{\infty} (-t) \, e^{izt} \, dt = \frac{1}{z^2}, 
\qquad 
\int_{0}^{\infty} \frac{-t^2}{2} \, e^{izt} \, dt = \frac{i}{z^3}\quad \text{if}~\Im(z)>0. 
\end{equation}
Hence, 
\begin{equation} \label{eq_0913_1}
\int_{0}^{\infty} g_F(t) e^{izt} \, dt 
= \frac{1}{z^2} \frac{\xi_F^\prime}{\xi_F}\left(\frac{1}{2}-iz \right)
\end{equation}
holds when $\Im(z)>1/2$. 

By 
\eqref{s202} and \eqref{s204}, 
the mapping $\gamma \mapsto -\overline{\gamma}$ 
defines a bijection from $\Gamma_F$ to $\Gamma_{F^\ast}$ with $m_\gamma=m_{-\bar{\gamma}}$, 
and $B_F=-B_{F^\ast}\,(\in \R)$. 
Applying these facts to \eqref{eq_101}, 
\begin{equation} \label{eq_0914_2}
\aligned 
\overline{g_F(t)} 
&= iB_F\,t - \frac{m_0}{2}\,t^2 + \sum_{0\not=\gamma\in \Gamma_F} m_\gamma
\frac{e^{i\overline{\gamma} t}-1}{\overline{\gamma}^2} \\
 &= -iB_{F^\ast}\,t - \frac{m_0}{2}\,t^2 
+ \sum_{0\not=\gamma\in \Gamma_{F^\ast}} m_\gamma \frac{e^{-i\gamma t}-1}{\gamma^2}
= g_{F^\ast}(t). 
\endaligned 
\end{equation}
Therefore, we obtain   
\begin{equation} \label{eq_0913_2}
\aligned 
2\int_{0}^{\infty} & \Re(g_F(t)) e^{izt} \, dt 
= \int_{0}^{\infty} (g_F(t)+g_{F^\ast}(t)) e^{izt} \, dt \\
& = \frac{1}{z^2} 
\left[ 
\frac{\xi_F^\prime}{\xi_F}\left(\frac{1}{2}-iz \right)
+
\frac{\xi_{F^\ast}^\prime}{\xi_{F^\ast}}\left(\frac{1}{2}-iz \right)
\right] 
 = \frac{1}{z^2} 
\frac{\xi_{FF^\ast}^\prime}{\xi_{FF^\ast}}\left(\frac{1}{2}-iz \right)
\endaligned 
\end{equation}
when $\Im(z)>1/2$.  

Let $\sigma_0:=\max\{\sigma \in \R\,|\, \xi_F(\sigma)=0\}$. 
This value is the same for $F^\ast$ by the functional equation \eqref{s202}. 
If $ \Re(-g_F(t))$ is non-negative for $t \geq t_0$, 
then the integral on the left-hand side converges 
when $\Im(z)>\max(\sigma_0-1/2,\,0)$ 
by \cite[Theorem 5b of Chapter II]{Wi41}. 
Hence, $\xi_F(s)$ and $\xi_{F^\ast}(s)$ have no zeros 
in the right-half plane 
$\Re(s)>\max(\sigma_0,\,1/2)$. 
In particular, the GRH for both $F$ and $F^\ast$ holds  
if $\xi_F(s)$ or $\xi_{F^\ast}(s)$ has no real zeros 
except for the possible zero at $s=1/2$. 
\hfill $\Box$

\subsection{Proof of Corollary \ref{cor_1}} 

If $F \in \mathcal{S}$, it has no zeros in $\Re(s)>1$ by the Euler product (S5). 
Further, there exists $\delta>0$ such that $\xi_F(s) \not =0$ for every $s \in (1-\delta,\infty)$ 
by the assumption. 
Therefore, the proof of Theorem \ref{thm_101} shows that 
$\xi_F(s)$ has no zeros in the right-half plane $\Re(s)>1-\delta$ 
due to the non-negativity of $\Re(-g_F(t)) \geq 0$ for $t \geq t_0$. 
Hence, $F(s) \not=0$ in the right-half plane $\Re(s)>1-\delta$. 
\hfill $\Box$ 
\medskip

This proof does not work for 
$F \in \mathcal{S}^\sharp \setminus \mathcal{S}$, 
since it may have a zero with real part greater than one.

\section{Zero-free formulas of screw functions} \label{section_4}

For functions $F$ in the semi-extended Selberg class $\mathcal{S}^{\sharp\flat}$, 
we obtain the following zero-free formulas of corresponding screw functions 
using the Lerch transcendent
\begin{equation} \label{eq_0123_2}
\Phi(w,s,a) := \sum_{n=0}^{\infty} \frac{w^n}{(n+a)^s}, 
\end{equation} 
where $a>0$, $|w|<1$, and $s \in \C$, or $|w|=1$ and $\Re(s)>1$. 


\begin{theorem} \label{thm_401}
Let $F \in \mathcal{S}^{\sharp\flat}$, and let $g_F(t)$ be the screw function defined by \eqref{eq_101}. 
We choose the data $\omega$, $Q$, $r$, $\lambda_j$, and $\mu_j$ 
satisfying (S3) for $F$ and 
define 
\begin{equation} \label{eq_401} 
\aligned 
\Psi_F(t)
& := 4\,\mathfrak{m}_F\,(e^{t/2}+e^{-t/2}-2) \\ 
& \quad  - \sum_{n \leq e^t} \frac{c_F(n)}{\sqrt{n}}(t-\log n) \\
& \quad 
+ \left[ \log Q 
+ \sum_{j=1}^{r} \lambda_{j} 
\frac{\Gamma'}{\Gamma}\left(\frac{\lambda_{j}}{2} +\mu_{j}\right) \right] t \\
&\quad  + 
\sum_{j=1}^{r} 
\lambda_{j}^2 \left[
\Phi\left(1,2,\frac{\lambda_{j}}{2}+\mu_{j}\right)
- e^{-t(\frac{1}{2}+\frac{\mu_{j}}{\lambda_{j}})}
\Phi\left(e^{-\frac{t}{\lambda_{j}}},2,\frac{\lambda_{j}}{2}+\mu_{j}\right)
\right]
\endaligned 
\end{equation}
for non-negative $t$ and set $\Psi_F(t):=\overline{\Psi_F(-t)}$ for negative $t$, 
where $c_F(n)$ are numbers in (S5'). 
Then $\Psi_F(t)$ is uniquely determined from $F$ and  
$\Psi_F(t)=-g_F(t)$ holds for all real numbers $t$. 
\end{theorem}
\begin{proof} 

To prove the equality $\Psi_F(t)=-g_F(t)$, it suffices to show that 
\begin{equation} \label{eq_0913_3}
-\int_{0}^{\infty} \Psi_F(t) e^{izt} \, dt 
= \frac{1}{z^2} \frac{\xi_F^\prime}{\xi_F}\left(\frac{1}{2}-iz\right) 
\end{equation}
holds when $\Im(z)>1/2$ by \eqref{eq_0913_1} and the uniqueness of the inverse Fourier transform. 

Taking the logarithmic derivative of \eqref{s203}, 
\[
\frac{\xi_F^\prime}{\xi_F}(s) 
= \frac{\mathfrak{m}_F}{s-1} + \frac{\mathfrak{m}_F}{s} + \frac{F'}{F}(s) 
+ \log Q + \sum_{j=1}^{r} \lambda_j \frac{\Gamma'}{\Gamma}(\lambda_j s +\mu_j). 
\] 
On the other hand, we have 
\begin{equation} \label{eq_0913_4} 
\int_{0}^{\infty} 4(e^{t/2}+e^{-t/2}-2)\, e^{izt} \, dt 
= -\frac{1}{z^2}\,\left(\frac{1}{(1/2-iz)-1}+\frac{1}{1/2-iz}\right) \quad \text{if}~\Im(z)>1/2, 
\end{equation}
\begin{equation*} 
\int_{0}^{\infty} \frac{(t-\log n)}{\sqrt{n}} \, \mathbf{1}_{[\log n,\infty)}(t)\, e^{izt} \, dt 
= -\frac{1}{z^2} \, n^{-(1/2-iz)} \quad \text{if}~\Im(z)>0, 
\end{equation*}
and \eqref{eq_0913_5} by direct and simple calculation. The second equality leads to 
\begin{equation} \label{eq_0913_6} 
\int_{0}^{\infty} 
\left(\sum_{n \leq e^t} \frac{c_F(n)}{\sqrt{n}}(t-\log n)\right) e^{izt} \, dt 
= \frac{1}{z^2}\left(-\sum_{n=2}^{\infty} \frac{c_F(n)}{n^{1/2-iz}} \right)
= \frac{1}{z^2}\frac{F'}{F}\left(\frac{1}{2}-iz\right)
\end{equation} 
if $\Im(z)>1/2$ by (S5'). 
Hence, \eqref{eq_0913_3} is proved if it is shown that 
\begin{equation} \label{eq_0913_7}
\aligned 
-\int_{0}^{\infty} & 
\lambda \left[
\Phi\left(1,2,\frac{\lambda}{2}+\mu\right)
- e^{-t\left(\frac{1}{2}+\frac{\mu}{\lambda}\right)}
\Phi\left(e^{-\frac{t}{\lambda}},2,\frac{\lambda}{2}+\mu\right)
\right] e^{izt} \, dt \\
&=
 \frac{1}{z^2} \left[
\frac{\Gamma'}{\Gamma}\left(\lambda \left(\frac{1}{2}-iz \right) +\mu\right)
-
\frac{\Gamma'}{\Gamma}\left(\frac{\lambda}{2} +\mu\right)
\right]
\endaligned 
\end{equation}
holds for $\lambda>0$, $\mu\in \C$ with $\Re(\mu)\geq 0$, and $\Im(z)>0$. 
We have 
\begin{equation*} 
\int_{0}^{\infty} \frac{1-e^{-\beta t}}{\beta^2} \, e^{izt} \, dt
=
\frac{i}{z^2}\left(\frac{1}{z+i\beta} -\frac{1}{i\beta} \right), 
\quad \Im(z) > \max(0,-\Re(\beta))
\end{equation*}
for a complex number  $\beta$ by direct and simple calculation. 
Therefore, 
\[
\aligned 
\int_{0}^{\infty} \frac{1-e^{-\frac{1}{\lambda}(\frac{\lambda}{2}+\mu+n) t}}
{\frac{1}{\lambda}(\frac{\lambda}{2}+\mu+n)^2} \, e^{izt} \, dt 
& =
\frac{1}{z^2}\left(\frac{1}{\lambda (1/2-iz) + \mu +n} -\frac{1}{n+1} \right) \\
& \qquad -\frac{1}{z^2}\left(\frac{1}{\frac{\lambda}{2}+\mu+n} -\frac{1}{n+1} \right). 
\endaligned 
\]
Using this and the series expansion \eqref{eq_0123_2}, 
the left-hand side of \eqref{eq_0913_7} is calculated as 
\begin{equation*} 
\aligned 
 = & -
\int_{0}^{\infty} 
\left[
\sum_{n=0}^{\infty}
\frac{1-e^{-\frac{1}{\lambda}(\frac{\lambda}{2}+\mu+n) t}}
{\frac{1}{\lambda}(\frac{\lambda}{2}+\mu+n)^2}
\right]\, e^{izt}\, dt  \\
& =- \frac{1}{z^2}\sum_{n=0}^{\infty}\left(\frac{1}{\lambda (1/2-iz) + \mu +n} -\frac{1}{n+1} \right) 
 +\frac{1}{z^2}\sum_{n=0}^{\infty}\left(\frac{1}{\frac{\lambda}{2}+\mu+n} -\frac{1}{n+1} \right).
\endaligned 
\end{equation*}
The right-hand side is equal to the right-hand side of \eqref{eq_0913_7} 
by the well-known series expansion 
\begin{equation*} 
\frac{\Gamma'}{\Gamma}(w) = -C_0 - \sum_{n=0}^{\infty}
\left( \frac{1}{w+n} - \frac{1}{n+1} \right), 
\end{equation*}
where $C_0$ is the Euler--Mascheroni constant. 
Hence, we complete the proof of \eqref{eq_0913_3}. 

On the other hand, it is known that 
the $\gamma$-factor $\gamma(s)=Q^s\prod_{j=1}^{r}\Gamma(\lambda_j s+\mu_j)$ 
is uniquely determined up to a constant multiple (\cite[Theorem 4.1]{Pe05}).  
Therefore, the logarithmic derivative $(\gamma'/\gamma)(s)$ 
is uniquely determined from $F$. 
Hence the values of $\Psi_F(t)$ are uniquely determined from $F$ by \eqref{eq_0913_3}. 
\end{proof}

\section{Complements} \label{section_5}

Finally, we make several comments on the generalizations of other results in \cite{Su22}. 

\subsection{} 

We have  
\begin{equation*} 
i\frac{\xi_F^\prime}{\xi_F}\left( \frac{1}{2}-iz\right) 
= - \left(B_F + \sum_{\gamma \in \Gamma_F\setminus\{0\}}\frac{m_\gamma}{\gamma} \right) 
- \frac{m_0}{z} 
- \sum_{\gamma \in \Gamma_F\setminus\{0\}}  \frac{m_\gamma}{z-\gamma} 
\end{equation*}
by \eqref{eq_0914_1}, 
where the sum is taken over $|\gamma|\leq T$ and then letting $T \to \infty$ 
for the convergence. By the functional equation \eqref{s202}, 
$\Gamma_F$ is closed under complex conjugation, 
thus $\sum_\gamma m_\gamma \gamma^{-1}$ is real. 
Therefore, the GRH for $F$ is true if and only if  
\[
\Im[\,i (\xi_F^\prime/\xi_F)(1/2-iz)\,] >0 \quad  \text{whenever $\Im(z)>0$}, 
\] 
because $\Im(-(z-\gamma)^{-1})=(\Im(z)-\Im(\gamma))/|z-\gamma|^2>0$ 
if $\Im(z)>\Im(\gamma)$ for each term 
(cf. \cite[Proof of Theorem 1.1]{La99} and \cite[Proof of Theorem 1]{La06}).  
Moreover, $y^{-1}(\xi_F^\prime/\xi_F)(1/2+y) \ll y^{-1}\log y$ 
as $y \to +\infty$ by (S5') and the Stirling formula of $(\Gamma'/\Gamma)(w)$. 
Combining these two facts,  
we obtain the following equivalence 
using the same argument as in the proof of \cite[Theorem 1.2]{Su22}. 

\begin{theorem} \label{thm_501}
Let $F$ be a member of 
$\mathcal{S}^{\sharp}$ 
and let $g_F(t)$ be the function defined by \eqref{eq_101}. 
Then, the GRH for $F$ is true if and only if 
$g_F(t)$ is a screw function on $\R$ in the sense of \cite{KrLa14}, 
that is, $g_F(-t)=\overline{g_F(t)}$ for all $t \in \R$ 
and $g_F(t-u)-g_F(t)-g_F(-u)+g_F(0)$ is a non-negative definite kernel on $\R$. 
\end{theorem} 

If $\gamma$ is a zero of $F(1/2-iz)$ for $F \in \mathcal{S}^\sharp$, 
its complex conjugate $\overline{\gamma}$ is also a zero with the same multiplicity, 
as follows from the functional equation \eqref{s202}. 
In addition, $B_F$ is real as observed in Section \ref{section_2_2}. 
Therefore, we have $g_F(t)=\overline{g_F(-t)}$ unconditionally. 
Thus, the problem is the non-negativity of the kernel $g_F(t-u) - g_F(t) - g_F(-u) + g_F(0)$. 
It is nothing but the non-negativity of the Weil distribution associated with $F$ 
as well as \cite[Section 3.4]{Su22}, and hence it is essentially nothing new. 
However, the equivalence in Theorem \ref{thm_501} helps us understand that 
Theorem \ref{thm_101} states the non-negativity of 
 $g_F(t-u) - g_F(t) - g_F(-u) + g_F(0)$ is reduced to the non-negativity of $\Re(-g_F(t))$ and the absence of real zeros of $F$.
\medskip

Assume the GRH for $F \in \mathcal{S}^{\sharp}$ and 
define the non-negative measure $\tau_F$ on the real line by 
$\tau_F  = \sum_{\gamma \in \Gamma_F} m_\gamma \delta_{-\gamma}$, 
where $\delta_a$ is the Dirac measure at the point $a \in \R$. 
Then $\int_{-\infty}^{\infty}(1+\lambda^2)^{-1}d\tau_F(\lambda)<\infty$ 
and $g_F(t)$ has the representation 
\[
g_F(t) = i\, C_F \,t + \int_{-\infty}^{\infty} 
\left( e^{i\lambda t} -1 - \frac{i\lambda t}{1+\lambda^2}  \right)
\frac{d\tau_F(\lambda)}{\lambda^2}
\]
on $\R$ with $C_F:=B_F+\sum_{0\not=\gamma\in \Gamma_F}m_\gamma \gamma^{-1}(1+\gamma^2)^{-1}$. 
This is nothing but the representation \cite[Theorem 5.1]{KrLa14} of screw functions on $\R$.

\subsection{} 

The theory of screw functions is somewhat simpler for real-valued function (\cite{KrLa14}).
The screw function $g_F$ for $F \in \mathcal{S}^\sharp$ 
is real-valued if $F=F^\ast$ by \eqref{eq_0914_2}. 
In general $F \in \mathcal{S}^\sharp$ does not satisfy $F=F^\ast$, 
but since $\mathcal{S}^\sharp$ is a multiplicative monoid, 
$G:=FF^\ast$ belongs to $\mathcal{S}^\sharp$ and satisfies $G=G^\ast$. 
Since the GRH for $F$ is equivalent to that for the product $FF^\ast$ 
by \eqref{s202}, it is sufficient to consider $FF^\ast$ to study the GRH. 
In other words, we only need to study $F \in \mathcal{S}^\sharp$ satisfying $F=F^\ast$.

\subsection{} 

The generalization of \cite[Theorem 1.6]{Su22} to $F \in \mathcal{S}^\sharp$ 
follows immediately from the integral representation \eqref{eq_0913_1}, 
but generalizing \cite[Theorem 1.7]{Su22} requires an additional assumption 
to prove the convergence of moments 
$\mu_{F,n}= \int_{0}^{\infty} e^{-t/2}(-g_F(t)) \, t^n \, dt$ 
for all non-negative integers $n$. 
If we suppose that $F \in \mathcal{S}$ and it has a zero-free region 
$\Re(s)>1-C(\log(3+|\Im(s)|))^{-1}$ for some $C>0$, 
then we obtain $(F'/F)(\sigma+it) \ll \log (3+|t|)$ 
in the zero-free region by standard arguments as in \cite[Section 5.6]{IK04}. 
With this estimate, 
the argument in \cite[Section 2.3]{Su22} works and 
the convergence of the moments $\mu_{F,n}$ is obtained.

\subsection{} 

In case $F \in \mathcal{S}^{\sharp\flat}$, 
the GRH for $F$ is true if and only if 
\[
\sum_{\gamma \in \Gamma_F}
\int_{-\infty}^{\infty} \phi(x) \,e^{i\gamma x} \, dx
\int_{-\infty}^{\infty} \overline{\phi(x)} \,e^{-i\gamma x} \, dx \geq 0
\]
for all smooth and compactly supported functions $\phi$ on $\R$ 
by Smajlovi\'{c}'s explicit formula \cite[Theorem 3.1]{Sm10}, 
the generalization of Li's criterion \cite[Theorem 4.3]{Sm10}, 
and argument of the proof of \cite[Theorem 1]{Bo01}. 
Therefore, 
the generalization of \cite[Theorem 1.3]{Su22} for $F \in \mathcal{S}^{\sharp\flat}$ holds 
by replacing $g$ with $g_F$. 
The result \cite[Theorem 1.4]{Su22} is also expected to be generalized 
to the semi-extended Selberg class $F \in \mathcal{S}^{\sharp\flat}$, 
but the author has not attempted it yet.
\medskip

\noindent
{\bf Acknowledgments}~
The author would like to thank the referee for carefully reading the original version and providing comments that were useful in improving the paper. 
This work was supported by JSPS KAKENHI Grant Number JP17K05163 and JP23K03050.

%

%

\bigskip 

\noindent
Masatoshi Suzuki,\\[5pt]
Department of Mathematics, \\
Institute of Science Tokyo \\
2-12-1 Ookayama, Meguro-ku, \\
Tokyo 152-8551, Japan  \\[2pt]
Email: {\tt msuzuki@math.titech.ac.jp}

\end{document}